\documentclass[11pt,
]
{article}

\usepackage{amssymb,amstext,amsthm,latexsym}
\usepackage[tbtags]{amsmath}
\usepackage{amsfonts}
\usepackage{bbm}

\usepackage{mathrsfs}
\usepackage{mathbbol}

\usepackage{mparhack}

\input{sol.sty}

\begin{document}

\title{A generalization of Solovay's $\Sigma$-construction\thanks
{This study was partially supported by RFBR grant 13-01-00006 
and Caltech.}}

\author{Vladimir Kanovei
}
%


\date{\today}
\maketitle


\begin{abstract}
A \ddd\Sg construction of Solovay is partially extended to 
the case of intermediate sets which are not necessarily subsets 
of the ground model.
As an application, we prove that, for 
a given name $t$, the set of all sets $t[G]$, 
$G$ being generic over the ground model, is Borel.
This result was first established by Zapletal by 
a totally different descriptive set theoretic argument.
\end{abstract}

\vspace{-3\baselineskip}
\def\contentsname{}
{\small\tableofcontents}


\punk{Introduction}

A famous \ddd\Sg construction by Solovay \cite{sol} shows
that if $x=t[G]$ is a real in a \ddd\dP forcing extension
$\mm[G]$ of a countable transitive model $\mm$, 
where $t\in\mm$ is a \ddd\dP name
and $\dP\in\mm$ is a forcing, then there is a set
$\Sg=\Sg(t,x)\in\mm[x]\yd \Sg\sq\dP$, such that
\ben
\renu
\itla{s1}
$G\sq\Sg$, 

\itla{s2}
if $G'\sq\Sg$ is \ddd\dP generic over $\mm$ then still $t[G']=x$, 

\itla{s3}
$\mm[G]$ is a \ddd\Sg generic extension of $\mm[x]$;
basically,
$G$ itself is \ddd\Sg generic over $\mm[x]$. 
\een
One might ask whether the \ddd\Sg construction can be
generalized to {\bf arbitrary} sets $x\in\mm[G]$, not
necessarily reals.
We don't know the answer, even in the easiest case 
$x\sq\pws\om$ not covered by Solovay's construction.

This note is devoted to a minor result in this direction.
Namely, given a countable transitive model $\mm$,
a forcing $\dP\in\mm$,
a \ddd\dP name $t\in\mm$,
and a set $X$ of any kind,
we define a set $\yha Xt$
(not a subset of $\dP$), which 
\ben
\tenu{\arabic{enumi})}
\item
satisfies a property which
resembles \ref{s1} (Lemma \ref{sgG}), \ and

\item is empty iff $X$ is not equal to 
$t[G]$ for any \ddd\dP generic set $G$ over $\mm$ 
(Corollary \ref{sf}).
\een 
As an application, we prove in the last section that, for 
a given \ddd\dP name $t$, the set of all sets $t[G]$, 
$G\sq\dP$ being generic over $\mm$, is Borel.
(Immediately, it is only analytic, of course.) 
This result was first established by Zapletal~\cite{t3} by 
a totally different argument.

\punk{Definition of $\mathbf\Sg$}

\bdf
[blanket assumptions]
\lam{bla}
We suppose that:\vom  

--- $\mm$ is a countable transitive model of $\ZFC$,\vom 

--- $\dP\in\mm$ is a forcing (a partial quasi-order),\vom

--- $t\in \mm$ is a \ddd\dP name of a transitive set
(so $\dP$ forces ``$t$ is transitive'').\vom

--- $X$ is a transitive set
(not necessarily $X\in \mm$ or $X\sq\mm$).\vom

\noindent
Let $\fo$ be $\fo_\dP^\mm$, the \ddd\dP forcing relation
over the ground model $\mm$. 
\edf

We also assume that a reasonable {\ubf ramified} system of names 
for elements of \ddd\dP generic extensions of $\mm$ is fixed. 
If $t$ is a name and $G\sq\dP$ is \ddd\dP generic over $\mm$ 
then let $t[G]$ be the \ddd Ginterpretation of $t$, so that 
we have 
$\mm[G]=\ens{t[G]}{t\text{ is a name}}$.
For any \ddd\dP names $s,t$, we let 
$s\prec t$ mean that $s$ occurs in $t$ as a name of a potential 
element of $t[G]$.
Then the set $\pe t=\ens{s}{s\prec t}$ 
(of all ``potential elements'' of $t$)
belongs to $\mm$ and
$$
t[G]=\ens{s[G]}{s\in\pe t\land \sus p\in G\,(p\fo s\in t)}\,.
$$
for any set $G\sq\dP$, generic over $\mm$.

If $d\sq\pe t$ then a condition $p\in\dP$
is called \rit{\ddd dcomplete\/} iff $p\fo s\in t$ holds for
all $s\in d$ and $p$ decides all formulas $s\in s'$ 
where $s,s'\in d$.

If $d$ is infinite then  \ddd dcomplete conditions do not 
necessarily exist.

\bdf
\lam{p*}
$\dpa(X,t)$ is the set of all pairs $\ang{p,a}$ such that 
$p\in\dP$, $u$ is a finite partial map, $\dom a\sq\pe t$, 
$\ran a\sq X$, and
$p$ is \ddd{(\dom a)}complete. 

We order $\dpa(X,t)$ so that $\ang{p,a}\le\ang{p',au'}$
($\ang{p,a}$ is stronger) 
iff $p\le p'$ in $\dP$ and $a$ extends $a'$ as a function.
\edf

Pairs in $\dpa(X,t)$ will be called \rit{superconditions}.
Given a supercondition  $\ang{p,a}\in\dpa(X,t)$, we'll call 
$p$ its \rit{condition}, and $a$ its \rit{assignment} --- 
because $a$ essentially assigns sets for (some) names which 
can be forced to be elements of $t$.
Note: generally 
speaking, superconditions are not members of $\mm$.

\bdf
\lam{defu}
Recall that $X$ is a fixed transitive set by Definition~\ref{bla}. 
Here we define a set $\yha Xt$ of all superconditions
$\ang{p,a}$ which, informally speaking, force nothing really
incompatible with the assumption that $X=t[G]$ for a set 
$G\sq\dP$ generic over $\mm$. 
The dependence on $\dP$ in the definition of $\yha Xt$ 
is suppressed.

\bit
\item
$\zha0Xt$ consists of all superconditions $\ang{p,a}\in\dpa(X,t)$ 
such that 

\quad 1) $\ran a\sq X$, and

\quad 2) if $s,s'\in\dom a$ then $p\fo s\in s'$ iff $a(s)\in a(s')$.

\item
If $\ga\in\Ord$ then the set 
$\zha{\ga+1}Xt$ consists of all superconditions
$\ang{p,a}\in\zha{\ga}Xt$ such that 

\quad
-- for any set $D\in M\yd D\sq\dP$, dense in $\dP$, 

\quad
-- and any name $s\in\pe t$, 

\quad
-- and any element $x\in X$, 

\noi
there is a stronger supercondition 
$\ang{q,b}\in \zha{\ga}Xt$ such that: 

\quad a) $\ang{q,b}\le\ang{p,a}$ and $q\in D$, 

\quad b) $x\in \ran b$, \ and

\quad c) \hspace{\mathsurround}either $s\in\dom b$ or $q\fo s\nin t$.

\item
Finally if $\la$ is a limit ordinal then 
$\zha{\la}Xt=\bigcap_{\ga<\la}\zha{\ga}Xt$.
\eit
The sequence of sets $\zha{\ga}Xt$ is decreasing,
so that there is an ordinal $\la=\la(X,t)$ such that
$\zha{\la+1}Xt=\zha{\la}Xt$; we let $\yha Xt=\zha{\la}Xt$.
\edf

The following is quite elementary.

\ble
\lam{sgu}
If\/ $\ang{p,a}\in\yha Xt$, a set\/ $D\in\mm\yd D\sq\dP$ is
dense in\/ $\dP$, and\/ $s\in\pe t$, $x\in X$,
then there is a pair\/ $\ang{q,b}\in \yha Xt$ satisfying$:$ 
$\ang{q,b}\le\ang{p,a}$, $q\in D$, $x\in \ran b$,
and either\/ $s\in\dom b$ or\/ $q\fo s\nin t$.
\ele
\bpf
This holds by definition, 
as $\yha Xt=\zha{\la}Xt=\zha{\la+1}Xt$.
\epf

We do not claim that if $\ang{p,a}\in\yha Xt$ and 
$q\in\dP\yt q\le p$ is a stronger condition then necessarily 
$\ang{q,a}\in\yha Xt$.
In fact this cannot be expected to be the case: 
indeed $q$ may strengthen $p$ in wrong way, that is, 
by forcing something 
incompatible with the assignment $a$.
Nevertheless, appropriate extensions of superconditions are 
always possible by Lemma \ref{sgu}.

\punk{Some results}

We continue to argue in the assumptions of
Definition~\ref{bla}.

\vyk{
\ble
\lam{uv}
If superconditions\/ $\ang{p,u}$ and\/ $\ang{p,v}$ 
(with the same\/ $p$) belong to\/ $\Sg(X)$ then\/ $u,v$ 
are compatible in the sense that\/ $u(s)=v(s)$ for all\/ 
$s\in\dom u\cap\dom v$.
\ele
\bpf
We argue by induction on the rank of $s$ in $\pe t$. 

If the rank is $0$ then $s$ is the name of the empty set. 
We claim that then $u(s)=v(s)=\pu$. 
Indeed suppose that $u(s)=x\ne\pu$ and let $y\in x$; then $y\in X$ 
by transitivity. 
By Lemma~\ref{sgu} there is a stronger supercondition 
$\ang{p',u'}\in\Sg(X)$ with $y\in\ran u'$, so that $y=u'(s')$, 
where $s'\in\dom{u'}\sq\pe t$. 
Thus $u'(s')\in u(s)=u'(s)$, so that by definition $q\fo s'\in s$, 
which contradicts the assumption that $s$ 
is the name of the empty  set.

Now the step. 
Suppose towards the contrary that $u(s)\ne v(s)$, 
say $y\in u(s)\bez v(s)$.
By Lemma~\ref{sgu} there is a pair of superconditions 
$\ang{p',u'}\yd\ang{p',v'}\in\Sg(X)$ (with the same $p'$) 
and names $s'\in \dom u'\yd t'\in\dom v'$ such that 
$u'(s')=v'(t')=y$. 
Moreover it follows from the hierarchical organization of names 
and still Lemma~\ref{sgu} 
\epf

}

Let $G\sq\dP$ be a \ddd\dP generic set over\/ $\mm$.
Say that a function (assignment) $a$, with $\dom a\sq\pe t$, is  
\rit{\ddd Gcompatible} if $a(s)=s[G]$ for all $s\in\dom a$.

The next lemma needs some work.

\ble
\lam{sgG}
Let\/ $G\sq\dP$ be a\/ \ddd\dP generic set
over\/ $\mm$, and\/ $t[G]=X$.
If\/ $\ang{p,a}\in\dpa(X,t)$, $a$ is\/ \ddd Gcompatible, 
and\/ $p\in G$,  then\/ $\ang{p,a}\in\yha Xt$.
\ele
\bpf
Prove $\ang{p,a}\in\zha\ga Xt$ by induction on $\ga$.

Assume that $\ga=0$.
By the \ddd{(\dom a)}completeness, 
if $s,s'\in \dom a$ then $p$ decides $s\in s'$.
If $p\fo s\in s'$ then $s[G]\in s'[G]$,
therefore $a(s)\in a(s')$ by the \ddd Gcompatibility.
Similarly, if $p\fo s\nin s'$ then $a(s)\nin a(s')$.

The step $\ga\to\ga+1$.
Suppose, towards the contrary, that 
$\ang{p,a}\nin\zha{\ga+1}Xt$ but 
$p\in\zha\ga Xt$ by the inductive hypothesis.
By definition, there exist: a set $D\in\mm\yd D\sq\dP$, 
dense in $\dP$, and elements 
$s\in\pe t$, $x\in X$, such that no supercondition
$\ang{q,b}\in \zha\ga Xt$ satisfies all of 
\bce
$\ang{q,b}\le\ang{p,a}$, \
$q\in D$, \ $x\in \ran b$, \ 
and either $s\in\dom b$ or $q\fo s\nin t$.
\ece
By the genericity, there is a condition $q\in G\cap D$, $q\le p$.
As $t[G]=X$, there is a finite assignment 
$b:(\dom b\sq \pe t)\to X$ 
such that 
\begin{quote}
$a\sq b$, $x\in\ran b$, 
$r[G]\in t[G]$ and $b(r)=r[G]$ for every name $r\in\dom b$,
and either $s[G]\nin t[G]$ or $s\in\dom b$. 
\end{quote}
There is a stronger condition $q'\in G\cap D$ such that 
if in fact $s[G]\nin t[G]$ then $q'\fo s\nin t$, and even more, 
$q'$ is \ddd{(\dom b)}complete.
Then $\ang{q',b}\in\zha\ga Xt$ by the inductive hypothesis, 
a contradiction.

The limit step is obvious.
\epf

\ble
\lam{gins}
If\/ $\ang{p,a}\in\yha Xt$ then there is a set\/ $G\sq\dP$, 
\ddd\dP generic over\/ $\mm$, 
and such that\/ $p\in G$ and\/ $t[G]=X$. 
\ele
\bpf
As the model $\mm$ is countable, Lemma~\ref{sgu} allows to define a 
decreasing sequence of superconditions $\ang{p_n,a_n}\in\yha Xt$,
$$
\ang{p,u}=\ang{p_0,a_0}\ge\ang{p_1,a_1}\ge\ang{p_2,a_2}\ge\dots\,,
$$
such that the sequence $\sis{p_n}{n\in\om}$ intersects every 
set $D\in\mm\yd D\sq\dP$, dense in $\dP$ --- hence it naturally 
extends to a generic set $G=\ens{p\in\dP}{\sus n\,(p_n\le p)}$, 
and in addition, the union $\vpi=\bigcup_na_n:\dom\vpi\to X$ 
of all assignments $a_n$ satisfies: 
\ben
\aenu
\itla{gins1} 
$\ran\vpi=X$, 
$\dom\vpi\sq\pe t$, \ and 

\itla{gins2} 
for any $s\in\pe t$\,: 

\quad either $s\in\dom\vpi$ --- then $s[G]\nin t[G]$, 

\quad or $q\fo s\nin t$ for 
some $q\in G$ --- then $s[G]\nin t[G]$.
\een
Due to the transitivity of both sets 
$t[G]=\ens{s[G]}{s\in\dom\vpi}$ and $X=\ran\vpi$, 
to prove that $t[G]=X$, it suffices to check that 
$\vpi(s)\in \vpi(s')$ iff $s[G]\in s'[G]$, 
for all names $s,s'\in\dom\vpi$. 
By the construction of $\vpi$, there is an index $n$ 
such that $s,s'\in\dom a_n$. 
By definition, condition $p_n\in G$ is \ddd{(\dom a_n)}complete, 
so $p_n$ decides $s\in s'$. 

If $p_n\fo s\in s'$ then $s[G]\in s'[G]$, and on the other hand, 
as $\ang{p_n,a_n}\in\zha0Xt$, we have 
$\vpi(s)=a_n(s)\in a_n(s')=\vpi(s')$.

Similarly, if $q\fo s\nin s'$ then $s[G]\nin s'[G]$ and 
$\vpi(s)\nin \vpi(s')$.
\epf

The next lemma shows that the ordinals $\la(X,t)$ 
as in Definition~\ref{defu} are bounded in $\mm$ whenever 
$\yha Xt\ne\pu$.

\ble
\lam{bou}
There is an ordinal\/ $\laa(t)\in\mm$ such that\/ 
$\la(t[G],t)<\laa(t)$ for every set\/ $G\sq\dP$,\/ \ddd\dP generic 
over\/ $\mm$. 
\ele
\bpf
Assume that a set $G\sq\dP$ is \ddd\dP generic over  $\mm$. 
Then $X=t[G]\in\mm[G]$, and hence $\la(X,t)$ is an ordinal in 
$\mm$, and its value is forced, over $\mm$, by a condition in $G$.
\epf

\bcor
\lam{sf}
Let\/ $X$ be a transitive set. 
The following are equivalent\/$:$
\ben
\aenu
\itla{sf1}
there is a set\/ $G\sq\dP$, \ddd\dP generic over\/ $\mm$, 
such that\/ $t[G]=X\;;$


\itla{sf2}
$\yha Xt\ne\pu\;;$

\itla{sf3}
$\zha{\laa(t)}Xt=\zha{\laa(t)+1}Xt\ne\pu$.
\een
\ecor
\bpf
Use Lemmas~\ref{sgG}, \ref{gins}, \ref{bou}.
\epf

We finish with a question. 
Let $\yh Xt=\ens{p\in\dP}{\sus u\,(\ang{p,u}\in \yha Xt}$.
Is it true that if\/ $G\sq\yh Xt$ is a\/ \ddd\dP generic set
over\/ $\mm$ then\/ $X=t[G]$?

\vyk{
\ble
\lam{=x}
If\/ $G\sq\yg X$ is a\/ \ddd\dP generic set
over\/ $\mm$ then\/ $X=t[G]$.
\ele
\bpf
Because of Lemma~\ref{gins}, it suffices to prove that
$t[G]=t[G']$ for any two \ddd\dP generic sets $G,G'\sq\yg X$.
We check that, even more, $s[G]=s[G']$ for each name
$s\in\pe t\cup\ans t$.

We argue by induction on the name rank of $s$ in $\pe t\cup\ans t$. 
If the rank is $0$ then $s$ is the name of the empty set,
so obviously $s[G]=s[G']=\pu$.

Now the step. 
Suppose towards the contrary that $s[G]\ne s[G']$, 
say $y\in s[G]\bez s[G']$. 
By the construction of names, there is a name $s'\in\pe t$ of the 
rank strictly below the rank of $s$, such that $y=s'[G]$. 
Then $y=s'[G']$ too by the inductive hypothesis. 
\epf
}

\vyk{
\begin{prF}[Lemma~\ref{ml}]\rm
By Corollary~\ref{semifin} and Lemma~\ref{bou}, 
the structures\/ $\stk{t[G]}\in$ and\/ $\stk{|K|}K$ are isomorphic 
if and only if $\Sg(K)\ne\pu$ 
if and only if $\Sg_{\la^\ast}(K)=\Sg_{\la^\ast+1}(K)\ne\pu$, and 
the last requirement can be expressed in Borel way.  
\epf
}

\punk{The property of being generic-generated is Borel}

Another consequence of Lemma~\ref{bou} and other results above 
claims that, in the assumptions of Definition~\ref{bla}, the 
set of all sets of the form $t[G]$, $G\sq\dP$ being generic 
over $\mm$, is Borel in terms of an appropriate coding, of 
all (hereditarily countable) sets of this form, by reals. 
This result was first established by Zapletal~\cite{t3} by 
a totally different argument using advanced technique of 
descriptive set theory.

In order to avoid dealing with coding in general setting, we 
present this result only in the simplest nontrivial 
(= not directly covered by Solovay's original result) 
case when sets $t[G]$ are 
(by necessity countable) sets of reals.

For a real $y\in\bn$, we let 
$R_y=\ens{(y)_n}{n\in\om}\bez \ans{(y)_0}$, where 
$(y)_n\in\bn$  and $(y)_n(k)=y(2^n(2k+1)-1)$ for all $n$ and $k$.
Thus $\ens{R_y}{y\in\bn}$ is the set of all at most 
countable sets $R\sq\bn$ (including the empty set).

\bte
\lam{zapt}
In the assumptions of Definition~\ref{bla}, if\/ $\dP$ forces 
that\/ $t[G]$ is a subset of\/ $\bn$ then the set\/ 
$W$ of all reals\/ $y\in\bn$, 
such that\/ $R_y=t[G]$ for a set\/ $G\sq\dP$ generic over\/ $\mm$, 
is Borel.
\ete
\bpf
Let $\vartheta$ be the least ordinal not in $\mm$.
By Corollary~\ref{sf}, for a real $y$ to belong to $W$ each of the 
two following conditions is necessary and sufficient:
\ben
\Renu
\itla{zapt1}
there exist an ordinal $\la<\vartheta$ and a sequence of sets 
$\zha\ga Xt\yt \ga\le{\la+1}$, where $X=R_y$, satisfying 
Definition~\ref{defu} and such that 
$\zha\la Xt=\zha{\la+1} Xt\ne \pu$;

\itla{zapt2}
for any ordinal $\la<\vartheta$ and any sequence of sets 
$\zha\ga Xt\yt \ga\le{\la+1}$, where $X=R_y$, satisfying 
Definition~\ref{defu}, if 
$\zha\la Xt=\zha{\la+1} Xt$ then 
$\zha\la Xt\ne \pu$.
\een
Condition \ref{zapt1} provides a $\fs11$ definition of 
the set $W$ 
while condition \ref{zapt2} provides a $\fp11$ definition 
of $W,$  
both relative to a real parameter coding the \ddd\in structure 
of $\mm$.
\epf

\end{document}